\documentclass[12pt]{amsart}
\usepackage{amssymb}
\usepackage{amsfonts}
\usepackage{amsmath}
\usepackage{graphicx}
\usepackage{tikz}
\usepackage{hyperref}
\usepackage[all]{xy}
\hypersetup{colorlinks = true,   urlcolor = blue,   linkcolor = blue,   citecolor = red}

\newtheorem{theorem}{Theorem}

\newtheorem{lemma}[theorem]{Lemma}
\newtheorem{notation}[theorem]{Notation}

\newtheorem{proposition}[theorem]{Proposition}
\newtheorem{remark}[theorem]{Remark}

\newtheorem{assumption}[theorem]{Assumption}

\usepackage{soul}
\usepackage[foot]{amsaddr}

\usepackage{wasysym}
\usepackage{xcolor}
\definecolor{lime}{HTML}{A6CE39}
\DeclareRobustCommand{\orcidicon}{%
	\begin{tikzpicture}
	\draw[lime, fill=lime] (0,0) 
	circle [radius=0.16] 
	node[white] {{\fontfamily{qag}\selectfont \tiny ID}};
	\draw[white, fill=white] (-0.0625,0.095) 
	circle [radius=0.007];
	\end{tikzpicture}
	\hspace{-2mm}
}
\foreach \x in {A, ..., Z}{%
	\expandafter\xdef\csname orcid\x\endcsname{\noexpand\href{https://orcid.org/\csname orcidauthor\x\endcsname}{\noexpand\orcidicon}}
}

\title[Explicit Homology Representations]{Explicit Homology Representation for Finite Groups Acting on Riemann Surfaces}
\date{June 2026}

\author[S. Allen Broughton]{S. Allen Broughton$^1$ \orcidB{}}
\address{
$^1$ Department of Mathematics \\ Rose-Hulman Institute of Technology \\ Terre Haute, Indiana, U.S.A.
}
\email{brought@rose-hulman.edu}

\author[Linden Disney-Hogg]{Linden Disney-Hogg$^{2,\ast}$ \orcidA{}}
\address{
$^2$ Department of Mathematics \\ Imperial College London \\ London, U.K.
}
\email[Corresponding author]{a.disneyhogg@imperial.ac.uk}

\thanks{$^\ast$ Corresponding author}

\thanks{\textbf{Acknowledgements.} LDH is grateful to Josh Fogg for their consultation on the efficiency of sparse matrix methods.}

\thanks{\textbf{Data Availability.} The datasets generated during the current study and the code for their creation/analysis are openly available at \url{https://github.com/DisneyHogg/homology_representation/}.}

\thanks{\textbf{Statements and Declarations.} The authors have no relevant financial or non-financial interests to disclose.}

\begin{document}
 
\begin{abstract}
Given a finite group $G$ acting orientably on a surface $S$ of genus $\sigma
\geq 2$, the group $G$ acts faithfully on the homology group $H_{1}(S;\mathbb{Z})$,
preserving the symplectic intersection form. The action on $S$ and the homology 
is determined by a \emph{generating vector}, a tuple of elements of $G$, 
generating $G$ and satisfying certain properties. In this note we show how 
to compute the homology representation, using the generating vector, 
when $S/G$ has genus 0 and the genus is suitably low. A $2\sigma \times 2\sigma$ 
representing matrix can be determined for any element in the group, usually
for a small set of generators. The matrices are computed with respect to 
an auto-generated basis for the cellular homology of $S$, using a regular 
$CW$ structure on $S$, derived from the $G$ action. 
% In addition, the matrix of intersection products of the auto-generated basis vectors is determined.
We demonstrate the application of these results by computing \emph{invariant theta characteristics} of the Riemann surfaces $S$ with the algorithm implemented using Sage. 
\end{abstract}

\maketitle

% \tableofcontents

\section{Introduction}\label{sec-intro}

\subsection{Problem description}
If a finite group $G$ acts orientably on a surface $S$, then $G$ acts on the
homology group $H_{1}(S;\mathbb{Z})$, preserving the symplectic intersection
form. If the surface is hyperbolic, i.e., of genus $\sigma \geq 2$, then the
action is faithful. For various applications, it is desirable 
to know this representation explicitly.
%, see Sections \ref{subsec-apps} and \ref{sec-theta}. 
In this note, we present an algorithm,
implemented in Sage \cite{Sa} and available from \url{https://github.com/DisneyHogg/homology_representation/}, to explicitly determine 
the representation by using simple geometric information about the quotient
surface $T=S/G$, and a generating vector for the action of $G$ (see Section \ref{subsec-action}). Our computational process auto-generates
a homology basis $a_1,\dots a_{2\sigma}$ adapted to the $G$ action. We are then
able to specify a $2\sigma \times 2\sigma $ representing matrix for any element
of the group, e.g., a small set of group generators.
% In addition, we show how to compute the intersection matrix using cup products which allows us to study the $G$-invariant symplectic form.
Our algorithm allows for homology coefficients from an arbitrary unital ring; the case of coefficients mod 2 is of special interest.

\subsection{Previous work}\label{subssec-prev-work}

There are existing methods for computing the homology representation from the data of the generating vector. Initially \cite{Gi1, Gi2} focused on the relatively simple case when $G$ is cyclic of prime order, later extending to arbitrary $G$ \cite{Gi3}, and worked with the fundamental polygon to find a basis of cycles adapted to the automorphism in order to simplify the expression of its action. Finding the adapted homology basis required the application of a rewriting system and repeated Teitze transformations. It was also shown how to compute the intersection form in this adapted basis \cite{GP}, which required finding a standard representative of each homology class. To our knowledge, no implementation of this algorithm was made publicly available. 

An alternative approach was developed in \cite{BeRoRo}, again using the fundamental polygon, now restricting to the case when the quotient surface has genus 0. There, the fundamental polygon for $S/G$ is lifted from the Riemann sphere to tessellate the surface by iteratively acting with the generators of the rotations of the fundamental polygon until the tessellation is complete. This tessellation is then treated as a simplicial complex for which the homology can be computed via standard methodologies. Again intersections of cycles are worked out by finding standard representatives of each homology class and carefully inspecting the order of intersections to determine the sign. Code was provided for this method, available at present from \url{https://github.com/rojas-ani/sage-routines}, and we use this for comparison. 

Neither of these two works discussed the additional problem of determining the action of $G$ on $H_1(S; R)$ when $R$ is a unital ring other than $\mathbb{Z}$. This may always be worked out from $H_1(S; \mathbb{Z})$ using the universal coefficient theorem, but the method we shall describe here makes clear how to use the ring $R$ from the beginning (providing computational benefit) while also avoiding some of the other difficulties inherent in former methods.

\subsection{Applications}\label{subsec-apps}

Explicit knowledge of the homology representation can be useful for a variety of reasons, including identifying loci in the moduli space of abelian varieties \cite[\S7]{BeRoRo}, and decomposing Jacobians. In this paper we shall consider the problem of computing the orbits of theta characteristics on $S$, which requires the homology representation with coefficients taken mod 2 \cite{Kallel2010, Braden2025}. This makes clear the benefit in having a method for computing the action which works over more general rings than $\mathbb{Z}$. In \S\ref{sec-theta} we shall demonstrate the performance of our algorithm in computing the homology representation mod 2 and see that it outperforms existing methodologies both in terms of runtime and memory usage. As such, we are able to compute the orbits of theta characteristics on curves which were previously unattainable, including certain modular curves. 

\subsection{Computation}\label{subsec-computation}
Our computations are done using Sage \cite{Sa}, though Magma \cite{Ma} or GAP \cite{Ga} would also work. The motivation for using Sage is its combination of being open-source, written using the clean language of Python, and retaining interfaces to many other mathematical software programs such as GAP, which itself was interfaced with for computation of the permutation action of $G$ on its cosets. 

\section{Preliminaries}\label{sec-prelim}
To formulate our problem precisely, we need to sketch the connection between
orientable group actions and the homology representation.
Throughout the paper, $S$ is a closed Riemann surface of genus $\sigma$, 
$G$ is a finite group that acts upon $S$ and $T=S/G$ is the quotient surface. 

\subsection{Orientable actions and signatures}\label{subsec-action}
The finite group $G$ \emph{acts orientably }on $S$ if there is a
monomorphism: 
\begin{equation}\label{eq-action}
\epsilon :G\rightarrow \mathrm{Homeo}^{+}(S),
\end{equation}
the group of orientation preserving homeomorphisms of $S$. If the image
consists of conformal automorphisms of $S,$ i.e., 
\begin{equation}\label{eq-Caction}
\epsilon :G\rightarrow \mathrm{Aut}(S),
\end{equation}
we say that $G$ \emph{acts conformally} on $S$. 
Every action as in \eqref{eq-action} gives rise to a representation 
on the homology $H_1(S;\mathbb{Z})$.
%Conformal actions \eqref{eq-Caction} also generate representations of the space of $q$-differentials $\mathcal{H}^q(S)$. The equivalence and characters of these various representations were studied in \cite{Bre, Br1, Br4} among other works.

For every orientable action there is a conformal structure on $S$, perhaps 
a multi-parameter family of such structures, preserved by $G$. These 
multi-parameter families are made up of so-called equisymmetric strata
in moduli space; the strata are connected, smooth quasi-projective varieties \cite{Br2}. 
It is well-known that we get the exactly one equivalence class of 
representations along a stratum \cite{Br4}.
It turns out that in each genus there are only a finite number of strata, 
a finite number of groups and a finite number of classes of representations to
consider \cite{Br4}. 

Two actions $\epsilon _{1},\epsilon _{2}$ of $G$ on possibly different surfaces $
S_{1},S_{2}$ are \emph{topologically equivalent} if there is an intertwining
homeomorphism $h:S_{1}\rightarrow S_{2}$ and an automorphism $\omega \in 
\mathrm{Aut}(G)$ such that:

\begin{equation}\label{eq-equivact}
\epsilon _{2}(g)=h\epsilon _{1}(\omega (g))h^{-1},\text{ }\forall g\in G,
\end{equation}
or in diagram form: 
\begin{equation}\label{dia-top-equiv} 
\xymatrix{
  G \ar[r]^-{\epsilon_2} \ar[d]^{\omega} & {\mathrm{Homeo}^+(S_2)}    \ar[d]^{Ad_h^{-1}} \\
  G \ar[r]^-{\epsilon_1}    &      {\mathrm{Homeo}^+(S_1)}
  }
\end{equation}
where $Ad_x(y) =xyx^{-1}$ whenever the composition is well defined. 
If $h$ is a conformal map then we say that the actions are conformally equivalent.

The quotient surface $T=S/G$ of an orientable action is a closed, orientable Riemann surface of genus $\tau $, and the quotient map 
\begin{equation}
\pi _{G}:S\rightarrow T \label{eq-Gquotient}
\end{equation}
is branched over $t$ \emph{branch points} $\{z_1,\ldots,z_t\} := B_G \subset T$. 
At each point $\tilde{z}_j$ lying over $z_j$, the stabilizer $G_{\tilde{z}_j}$ is a cyclic group 
of order $n_j$. We call $\mathfrak{s} =(\tau;n_1,\ldots,n_t)$ the \emph{signature}
of the action. The relation between the genus, the group order, and the signature
is given by the Riemann-Hurwitz equation
\begin{equation}\label{eq-RH}
\frac{2\sigma-2}{|G|} = 2\tau-2+t-\sum_{j=1}^t\frac{1}{n_j}.
\end{equation} 

\subsection{The homology representation}\label{subsec-homrep}
Given a group action $\epsilon$, we define the homology representation 
\begin{equation*}
\rho :G\rightarrow Sp_{2\sigma}(\mathbb{Z})
\end{equation*}
by setting $\rho (g)=\epsilon (g)_{\ast }$ to be the automorphism of 
$H_{1}(S;\mathbb{Z})$ induced by the homeomorphism $\epsilon (g)$ of $S$.
Since $G$ preserves the non-degenerate intersection form on $H_{1}(S;\mathbb{
Z}),$ and $H_{1}(S;\mathbb{Z})$ is a free $\mathbb{Z}$-module the
representation has its image in the integral symplectic group 
$Sp_{2\sigma}(\mathbb{Z})$.
In a suitable basis the matrices have the form%
\begin{equation*}
\rho (g)=\left[ 
\begin{array}{cc}
A & B \\ 
-B & A%
\end{array}
\right]
\end{equation*}
for matrices $A,B\in SL_{\sigma}(\mathbb{Z})$. If two actions 
$\epsilon _{1},\epsilon _{2}$ of $G$ are equivalent as in equation 
\eqref{eq-equivact} then the homology representations $\rho _{1},\rho _{2}$ are
equivalent in the sense that 
\begin{equation}
\rho _{2}(g)=U\rho _{1}(\omega (g))U^{-1},\text{ }\forall g\in G,
\label{eq-equivrep}
\end{equation}%
where the matrix $U$ is the the induced map $h_\ast: H_1(S_1;\mathbb{Z})\rightarrow H_1(S_2;\mathbb{Z})$ with respect to the corresponding bases.

\begin{remark}\label{rk-inner}
Note that if $\omega=Ad_g$ is an inner automorphism then it can be absorbed 
into the matrix $U$ i.e., $U^\prime = U\rho_1(g)$. Therefore, the equivalence 
class of a representation is not affected by inner automorphisms. 
Since topological equivalence involves the entire automorphism group of $G$ 
we do have to consider the $\mathrm{Out}(G)$ action upon the representations.
\end{remark} 

\begin{remark}\label{rk-coeffR}
Of course, $G$ acts upon the homology $H_{1}(S;R)$ and cohomology $H^{1}(S;R)$ 
when $R$ is any unital ring. The case of most importance to application is 
mod 2 coefficients where $R=\mathbb{Z}_2$, see Remark \ref{rk-gencoeff} and \S\ref{sec-theta}. 
\end{remark} 

\subsection{Fundamental groups and generating vectors}\label{subsec-Pi1GV}
\paragraph{\textbf{Actions and monodromies.}}
It is typical practice to use Fuchsian group pairs $\Pi \vartriangleleft \Gamma$
and surface kernel epimorphisms 
\begin{equation}
\Pi \hookrightarrow \Gamma \overset{\eta }{\twoheadrightarrow }G
\label{eq-uniform1}
\end{equation}%
to analyze $G$-actions.
However, since we are considering actions in the topological category, we will
utilize an equivalent method using fundamental groups and covering spaces for 
our constructions.

Set $T^{\circ }=T-B_{G}$ and let $S^{\circ }=\pi _{G}^{-1}\left( T^{\circ }\right)$.
Then the restricted map 
\begin{equation}\label{eq-unramifiedcover}
 \pi _{G}:S^{\circ }\rightarrow T^{\circ }
\end{equation}
is a regular, unramified covering space whose group of deck transformations 
\[
\mathrm{Gal}(\pi_G) = \mathrm{Gal}(S^\circ/T^\circ) = \{\phi \in \mathrm{Aut}(S^\circ) : \pi_G\circ\phi = \pi_G \}
\]
equals $\epsilon (G)$, restricted to $S^{\circ }$. This covering determines a normal
subgroup 
$\Pi_{G}=\pi_{1}(S^{\circ })\vartriangleleft \pi_{1}(T^{\circ })$
and an exact sequence 
\[
\Pi _{G}\hookrightarrow \pi _{1}(T^{\circ})\overset{\lambda}{\twoheadrightarrow} 
\mathrm{Gal}(S^\circ/T^\circ)=\epsilon(G) 
\]
by lifting loops to deck transformations. 
Combine the map $\lambda$ with $\epsilon (G)\overset{\epsilon ^{-1}}{\rightarrow }G$ 
to get an exact sequence 
\begin{equation}\label{eq-monodromyexactsequence} 
\Pi_{G}\hookrightarrow \pi_{1}(T^{\circ })\overset{\xi }{\twoheadrightarrow} G. 
\end{equation}
Any surjective map $\xi: \pi _{1}(T^{\circ })\rightarrow G$ is called a \emph{monodromy}.
Monodromies and generating vectors, discussed below, will be our basic starting point for specifying actions.

Since we have left out base points to simplify the exposition, $\lambda$ and $\xi$ are 
ambiguous up to inner automorphisms. However, this is no concern according to
Remark \ref{rk-inner}, since we only looking at homology actions up to inner automorphisms. 
The relation among $\epsilon$, $\lambda$ and $\xi$ is 
\begin{equation}\label{eq-eps-xi}
\xi=\epsilon^{-1}\circ\lambda \text{ or } \epsilon = \lambda\circ \xi^{-1}.  
\end{equation}
The right hand side of the equation makes sense as any two $\xi$-preimages $\gamma$, $\gamma^\prime$ of $g \in G$ satisfy 
$\gamma^\prime=\gamma\delta$ with $\delta \in \Pi_G$, and the elements of 
$\Pi_{G}=\pi_{1}(S^{\circ })$ have only trivial lifts.  
For a brief discussion on base points and path lifting, see Section \ref{subsec-lifting}.

In the opposite direction, suppose we start with an exact sequence as in equation 
\eqref{eq-monodromyexactsequence}. Selecting a conformal structure on $T$
and using the exact sequence, we can construct an unramified, holomorphic, regular 
covering space which we still denote $\pi _{G}:S^{\circ }\rightarrow T^{\circ }$.
The deck transformations in $\mathrm{Gal}(\pi_G)$, found by path lifting, 
are automatically holomorphic. We define $\epsilon: G \rightarrow \mathrm{Aut}(S^\circ)$
by the second equation in \eqref{eq-eps-xi}. We can fill in the punctures to get a closed surface and 
the restricted action $\epsilon$ extends to a conformal action 
$\epsilon :G\rightarrow \mathrm{Aut}(S)$ at the filled in punctures 
using the Removable Singularity Theorem.

\paragraph{\textbf{Generating vectors.}}
The fundamental group $\pi _{1}(T^{\circ })$ has the following presentation: 
\begin{eqnarray}
\mathrm{generators} &:&\{\alpha _{i},\beta _{i},\gamma _{j} \, | \, 1\leq i\leq \tau
,1\leq j\leq t\}, \label{eq-piTpres-1} \\
\mathrm{relations} &:&\prod_{i=1}^{\tau }[\alpha _{i},\beta
_{i}]\prod_{j=1}^{t}\gamma _{j}=1. \label{eq-piTpres-2}
\end{eqnarray}
Define 
\begin{equation*}
a_{i}=\xi (\alpha _{i}),b_{i}=\xi (\beta _{i}),c_{j}=\xi (\gamma _{j}),
\end{equation*}
then the $2\tau +t$ tuple 
$\left( a_{1},\ldots ,a_{\tau },b_{1},\ldots,b_{\tau },c_{1},\ldots ,c_{t}\right)$ 
is called a \emph{generating vector} for the action. We observe that
\begin{equation} \label{eq-GV1}
G=\left\langle a_{1},\ldots ,a_{\tau },b_{1},\ldots ,b_{\tau},
c_{1},\ldots ,c_{t}\right\rangle,
\end{equation}
\begin{equation} \label{eq-GV2}
o(c_{j})=n_{j}, \text{ and }
\end{equation}%
\begin{equation} \label{eq-GV3}
\prod_{i=1}^{\tau }[a_{i},b_{i}]\prod_{j=1}^{t}c_{j}=c_{1}^{n_{1}}=\cdots
=c_{t}^{n_{t}}=1 
\end{equation}
for some integers $n_{j}$ $\geq 2$. We call $(\tau ; n_{1},\ldots,n_{t})$
the signature of the action. This definition agrees with that given at the end of Section \ref{subsec-action}.  

\begin{proposition}\label{prop-xiGV equiv}
Once an ordered generating set for $\pi_1(T^\circ)$, as in \eqref{eq-piTpres-1}, has been fixed, there is a 1-1
correspondence between monodromies and generating vectors.
\end{proposition}

\subsection{Path lifting and dependence on base points}\label{subsec-lifting} 
We make a few remarks on path lifting and base points as they will be needed
for the computation of equivariant cellular homology in Section
\ref{sec-maphomology}. Given $\epsilon$, the maps $\lambda$ and $\xi$ 
are dependent on a choice of a base point $w_0 \in T^\circ$ and lift point 
$\widetilde{w_0}$ lying over $w_0$. Specifically, for 
$\gamma \in \pi_1(T^\circ,w_0)$, we define the fibre action on the 
$\pi_G^{-1}(w_0)$ by:
\begin{equation}\label{eq-lift}
 \lambda(\gamma)\cdot\widetilde{w_0} = \widetilde{\gamma}(1),
\end{equation}
where $\widetilde{\gamma} : [0,1] \to S^\circ$ is the $\pi_G$-lift of $\gamma$ satisfying 
$\widetilde{\gamma}(0)=\widetilde{w_0}$. Since the cover \eqref{eq-unramifiedcover}
is regular $\lambda(\gamma)$ is actually a deck transformation in 
$\mathrm{Gal}(\pi_G)=\epsilon(G)$. It can be shown that if another lift point 
$\widetilde{w_0}^\prime$ is selected then the new lifting action $\lambda^\prime$ 
satisfies $\lambda^\prime(\gamma) = \epsilon(x)\lambda(\gamma)\epsilon(x)^{-1}$ for 
any $\gamma \in \pi_1(T^\circ,w_0)$ and some $x \in G$. In fact we may chose 
$x$ to satisfy $\epsilon(x)= \lambda(\gamma^\prime)$ where the loop $\gamma^\prime$
is chosen so that $\widetilde{\gamma^\prime}$ is a path from 
$\widetilde{w_0}^\prime$ to $\widetilde{w_0}$. It follows from equation
\eqref{eq-eps-xi} that the choice of a lift point affects the monodromy $\xi$ only by an 
inner automorphism. If the monodromy $\xi$ is given first, then the action $\epsilon$
is modified by an inner automorphism of $\mathrm{Gal}(\pi_G)$. 

Now consider the effect of varying the base point $w_0$. Let $w_1 \in T^\circ$ be 
some other base point and $\delta$ a path in $T^\circ$ from $w_0$ to $w_1$. 
Then $\mathrm{Ad}_\delta : \gamma \rightarrow \delta*\gamma*\delta^{-1}$ is an 
isomorphism 
\[
\mathrm{Ad}_\delta: \pi_1(T^\circ,w_1) \rightarrow \pi_1(T^\circ,w_0).
\]
The lifting homomorphisms $\lambda_0, \lambda_1$, based at $\widetilde{w_0},\widetilde{w_1}$,
respectively, fit with the map $\mathrm{Ad}_\delta$ into the diagram:
\begin{equation}\label{dia-lift2Gal} 
\xymatrix{
  \pi_1(T^\circ,w_1) \ar[r]^-{\lambda_1} \ar[d]^{\mathrm{Ad}_\delta} & \mathrm{Gal}(S^\circ/T^\circ)    \ar[d]^{Ad_\phi} \\
  \pi_1(T^\circ,w_0) \ar[r]^-{\lambda_0}               & \mathrm{Gal}(S^\circ/T^\circ)
  }
\end{equation}
i.e., for some  $\phi \in \mathrm{Gal}(S^\circ/T^\circ)$ 
\begin{equation}\label{eq-baseptEQ}
\lambda_0(\delta*\gamma*\delta^{-1}) = \phi \lambda_1(\gamma)\phi^{-1}.
\end{equation}
 
To see this, first consider the case where $\tilde{\delta}$ is a path from $\widetilde{w_0}$ 
to $\widetilde{w_1}$. The lift $\widetilde{\delta*\gamma*\delta^{-1}}$ has three stages:
\begin{itemize}
 \item $\widetilde{\delta}$ from $\widetilde{w_0}$ to $\widetilde{w_1}$,
 \item $\widetilde{\gamma}$ from $\widetilde{w_1}$ to $\lambda_1(\gamma)\cdot \widetilde{w_1}$, and
 \item $(\widetilde{\delta^{-1}})^\prime$ from $\lambda_1(\gamma)\cdot \widetilde{w_1}$ to 
    $\lambda_0(\delta*\gamma*\delta^{-1})\cdot \widetilde{w_0}$,
\end{itemize}
where we use a prime to indicate that the lift does not start at the standard starting point. 
The two paths $\lambda_1(\gamma)\cdot\widetilde{\delta^{-1}}$ and $(\widetilde{\delta^{-1}})^\prime$
both project to $\delta^{-1}$ and have the same starting point. So, by unique path lifting and since
$ \mathrm{Gal}(S^\circ/T^\circ)$ acts without fixed points, then 
\[
\lambda_0(\delta*\gamma*\delta^{-1})=\lambda_1(\gamma).
\]
So, in this case, equation \eqref{eq-baseptEQ} holds with $\phi=\mathrm{Id}$. 

Now suppose that $\tilde{\delta}$ is a path from $\widetilde{w_0}$ 
to $\widetilde{w_1}^\prime \neq \widetilde{w_1}$. Construct a path 
$\delta^\prime \in \pi_1(T^\circ,w_1)$ with a lift $\widetilde{\delta^\prime}$ 
from $\widetilde{w_1}^\prime$ to $\widetilde{w_1}$. By the previous case
\[
\lambda_0((\delta*\delta^\prime)*\gamma^\prime*(\delta*\delta^\prime)^{-1}) = \lambda_1(\gamma^\prime),
\]
for any $\gamma^\prime \in \pi_1(T^\circ,w_1)$.
Pick $\gamma^\prime = (\delta^\prime)^{-1}*\gamma* \delta^\prime$ and we get
\[
\lambda_0(\delta*\gamma*\delta^{-1}) = \lambda_1((\delta^\prime)^{-1}*\gamma*\delta^\prime))
 =\lambda_1(\delta^\prime)^{-1}\lambda_1(\gamma)\lambda_1(\delta^\prime).
\]
We pick $\phi = \lambda_1(\delta^\prime)^{-1}$ to get equation \eqref{eq-baseptEQ}.

\section{Equivariant cell-complex homology on S}\label{sec-EQhom}
Having explained the theory of how a generating vector determines monodromy, we shall explain how we can use this to compute the homology of $S$ and the action of $G$ upon it in the framework of $CW$-complexes and cellular homology. For a reference on $CW$-complexes, computing cellular homology, and the cellular boundary formula see \cite{Ha}. The discussion shall follow the development in \cite{BrCoIz3}. 

\subsection{Setting up the equivariant tiling on S}\label{subsec-tilingsetup}
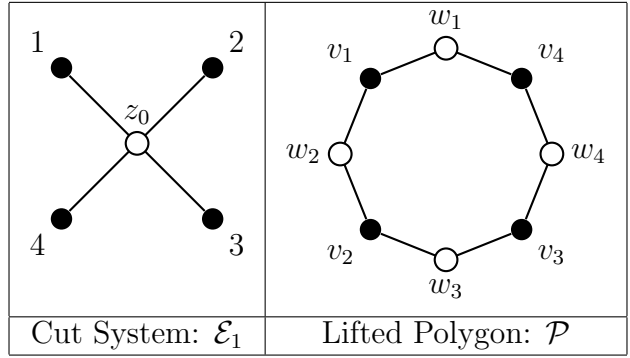
\begin{figure}[h!]
\begin{center}
\begin{tabular}
[c]{|c|c|}\hline

%%%% tikzpic E4
\begin{tikzpicture}[scale=0.5, inner sep=1mm]

\node (A) at (-2,-2) [shape=circle, fill=black] [label={-135:$4$}] {};
\node (B) at (-2,2) [shape=circle, fill=black] [label={135:$1$}] {};
\node (C) at (2,2) [shape=circle, fill=black] [label={45:$2$}] {};
\node (D) at (2,-2) [shape=circle, fill=black] [label={-45:$3$}] {};
\node (W) at (0,0) [shape=circle, fill=white] [label=$z_0$] {};
\node (P) at (0,-4) [shape=circle, fill=white] {}; %positioning

\draw [thick] (W) to (A);
\draw [thick] (W) to (B);
\draw [thick] (W) to (C);
\draw [thick] (W) to (D);
\filldraw[color=black, fill=white, thick] (W) circle (0.3);

\end{tikzpicture}

&
%%%% tikzpic P4
\begin{tikzpicture}[scale=0.5, inner sep=1mm]

\node (A) at (-2,-2) [shape=circle, fill=black] [label={-135:$v_2$}] {};
\node (B) at (-2,2) [shape=circle, fill=black] [label={135:$v_1$}] {};
\node (C) at (2,2)  [shape=circle, fill=black] [label={45:$v_4$}] {};
\node (D) at (2,-2) [shape=circle, fill=black] [label={-45:$v_3$}] {};

\node (E) at (-2.8,0) [shape=circle, fill=white] [label = {180:$w_2$}] {};
\node (F) at (0,2.8) [shape=circle, fill=white] [label = $w_1$] {};
\node (G) at (2.8,0) [shape=circle, fill=white] [label = {0:$w_4$}] {};
\node (H) at (0,-2.8) [shape=circle, fill=white] [label = {-90:$w_3$}] {};

\draw [thick] (A) to (E);
\draw [thick] (E) to (B);
\draw [thick] (B) to (F);
\draw [thick] (F) to (C);
\draw [thick] (C) to (G);
\draw [thick] (G) to (D);
\draw [thick] (D) to (H);
\draw [thick] (H) to (A);

\filldraw[color=black, fill=white, thick](-2.8,0) circle (0.3);
\filldraw[color=black, fill=white, thick](0,-2.8) circle (0.3);
\filldraw[color=black, fill=white, thick](2.8,0) circle (0.3);
\filldraw[color=black, fill=white, thick](0,2.8) circle (0.3);

\end{tikzpicture}
\\\hline
Cut System: $\mathcal{E}_1$ & Lifted Polygon: $\mathcal{P}$
\\\hline
\end{tabular}
\end{center}
\caption{Cut system and polygon}
\label{tikz-EnP}
\end{figure}

\paragraph{\textbf{Lifting maps on surfaces.}} 
We can create a $G$-equivariant cell complex on $S$, suitable for computing
homology, by lifting a \emph{map}, i.e. a $2$-dimensional cell complex, 
from $T$ to $S$. To this end, let $\mathcal{E}$ be a connected, undirected graph,
embedded in $T$, such that the following hold. 
\begin{enumerate}
 \item The \emph{vertex set} of $\mathcal{E}$ contains all the branch points,
  $B_G$, and possibly some regular points.
 \item The \emph{edges} of $\mathcal{E}$ are smooth loops or undirected arcs that meet at
  vertices only and the tangent directions at a vertex are all distinct.
 \item The complement $T-\mathcal{E}$ is a disjoint union of open topological disks 
  called \emph{faces}, using graph theoretic terminology.
  We think of the faces as open polygons whose ``ideal boundary'' consists of edges, concatenated to form a topological circle. 
  The actual boundary may be quite different.
\end{enumerate}
An example of such a graph is given in the left panel of Figure \ref{tikz-EnP}, where it is called 
a \emph{cut-system}. The complement has a single face. A model of the corresponding
open polygon (and its homeomorphic lift to $S$) is shown in the right panel, along 
with its ideal boundary. The ramified nodes are numbered and coloured black. 
The regular nodes are white. The labelling of the polygon vertices will be explained 
in Section \ref{subsec-1face}.
 
\begin{remark}\label{rk-map}
For our purposes, any surface-graph pair $(U,\mathcal{G})$ is a \emph{map},
provided that the following conditions are met. The surface $U$ is closed and
orientable and the graph $\mathcal{G}$ is a connected, undirected graph, 
embedded in $U$ satisfying items 2 and 3 above. We do not consider 
the case where $U$ is not orientable. 
\end{remark} 

\begin{notation}\label{not-tilde}
The lift to $S$ of a vertex, edge, face, or subgraph of $\mathcal{E}$ may be 
denoted by decorating the object with a tilde, as previously done in \S\ref{sec-prelim}. For example, $\widetilde{\mathcal{E}}$ 
is the lift of $\mathcal{E}$ to $S$ and $(S, \widetilde{\mathcal{E}})$ is the 
lift of the map $(T,\mathcal{E})$. To keep notation simple, we only use the tilde
decoration to eliminate confusion, as in the examples $\widetilde{\mathcal{E}}$ and
$(S, \widetilde{\mathcal{E}})$. 
\end{notation}

The lift $\widetilde{\mathcal{E}}$ of $\mathcal{E}$ to $S$ defines a map on $S$.
For, the cover $\pi_G$ introduced in equation \eqref{eq-unramifiedcover} is
unramified away from $B_G$ and the inverse image of an open face in $T-\mathcal{E}$ 
is a disjoint union of $|G|$ homeomorphic lifts of the given face to 
$S-\widetilde{\mathcal{E}}$, Likewise, each open edge in $\mathcal{E}$ 
lifts to $|G|$ open edges in $\widetilde{\mathcal{E}}$. Because the cover $\pi_G$ 
is regular with deck transformation group $\epsilon(G)$, $G$ acts freely on the sets
of open faces and open edges of $S$.

In \cite{BrCoIz3}, the case of lifting a single-faced map $(T,\mathcal{E})$ 
to $(S,\widetilde{\mathcal{E}})$ was studied in detail. The graph $\mathcal{E}$ 
was called a \emph{cut system} and the components of the lifted open face 
were called \emph{polygons} (again see Figure \ref{tikz-EnP}). 
Criteria for determining when the lifted map $(S,\widetilde{\mathcal{E}})$ 
was a regular $CW$-complex (see Section \ref{subsec-CWcomplex}) were specified. 
In the regular $CW$ case, for each open edge $e \in \mathcal{E}$, exactly 
two of the edges in the total lift $\pi_G^{-1}(e)$ are included in the 
boundary of a lifted polygon $\mathcal{P}$. Finally, if $T$ is topologically a sphere, 
then the cut system $\mathcal{E}$ is a tree.

%%%% tikzfigure cut-sytem and polygon

\subsection{CW-complexes}\label{subsec-CWcomplex} 
Let us set up our notation for the structure of a $CW$ complex $X$ so that 
we may use it in the calculation of cellular homology. 
A general $CW$-complex $X$ is given by the disjoint union: 
\[
X = \bigcup_{n}\bigcup_{i \in I_n}C^n_i, 
\]
where $C^n_i$ is an open cell of dimension $n$, homeomorphic to the interior 
of $B^n_i$, a copy of the closed unit ball in $\mathbb{R}^n$.
There is a separate index set $I_n$ for each dimension $n$. 
The $n$-skeleton is the union of cells of dimension $n$ or lower:
\[
X^n = \bigcup_{m\le n}\bigcup_{i \in I_m}C^m_i. 
\]
We denote the attaching map of $C^n_i$ by
\[ \alpha^{n}_i:\partial{B^n_i}\rightarrow X^{n-1}.
\]

In a \emph{regular $CW$-complex} $X$, we require that the closure in $X$ of each 
open $n$-dimensional cell $C^n_i$ be homeomorphic to the closed unit ball in 
$\mathbb{R}^n$, i.e., the embedding $C^n_i\rightarrow X^n$ extends to an embedding
of the closed ball $B^n_i\rightarrow X^n$. Though regularity is not needed 
for computing cellular homology and the $G$ action on homology, for our purposes 
it is convenient and no cost to assume regularity for computing boundary maps in 
cellular chain complexes. 

For the complex $X$ determined by the map $(U,\mathcal{G})$, $X^0= \text{vertex set}$,
$X^1=\mathcal{G}$, $X^2=U$ and all higher skeletons equal $U$. Regularity is easy 
to determine for a map $(U,\mathcal{G})$, given by the following condition. 
\begin{enumerate}
 \item Dimension 1: The closure of each $1$-cell in $\mathcal{G}$ is an undirected arc or loop
  with no internal self intersections. For regularity, there can be no loops.
 \item Dimension 2: For each closed face $f$ in $U$, the boundary, $\partial f$, is
   a closed subgraph of $\mathcal{G}$. For regularity, the boundary must be a closed
   cycle in $\mathcal{G}$, there can be no self intersections.
\end{enumerate}

\begin{remark}\label{rk-liftregularmap}
If $(T,\mathcal{E})$ is regular as a $CW$ complex then the lift 
$(S,\widetilde{\mathcal{E}})$ is also regular as a $CW$ complex.
\end{remark}

\subsection{Computing cellular homology}
We develop our process for computing equivariant homology for any map on a surface,
though when we get to computing specific examples, we will consider only planar 
actions (where $T$ has genus 0) and will choose lifts of specific maps on the sphere with one or two faces. 
This will simplify computation and directly relate the computation to the 
generating vector. 

There is a standard way to compute the homology groups of a $CW$-complex, $X$, which we will apply to the case of a map $(U,\mathcal{G})$. 
% See \cite{Ha} for background and details. 
We may use the \emph{cellular chain complex} 
\begin{equation}\label{dia-cellchain1} 
\xymatrix{
  \cdots H_{n+1}(X^{n+1},X^{n}) \ar[r]^-{\partial_{n+1}} & H_n(X^n,X^{n-1}) 
  \ar[r]^-{\partial_n} & H_{n-1}(X^{n-1},X^{n-2}) \cdots
  }
\end{equation}
to compute homology, via the isomorphism
\begin{equation}\label{eq-cellhomology1}
H_n(X) \simeq \mathop{ker}(\partial_n)/\mathop{im}(\partial_{n+1}).      
\end{equation} 
Of course, we need to identify $H_n(X^n,X^{n-1})$ and the differentials $\partial_n$.  
By using the homology excision theorem, we have: 
\begin{equation}\label{eq-cellhomology2}
H_n(X^n,X^{n-1}) \simeq \bigoplus_{i \in I_n}  H_n(B^n_i,\partial B^n_i)= 
\bigoplus_{i \in I_n} \mathbb{Z}\beta^n_i, 
\end{equation}  
where  $\beta^n_i$ represents a positively oriented generator of 
$H_n(B^n_i,\partial B^n_i)$. If we impose a geometric orientation 
on the $n$ cells the $\beta^n_i$ will be determined. 

Alternatively, 
\begin{equation}\label{eq-cellhomology3}
 H_n(X^n,X^{n-1}) =  \tilde{H}_n(X^n/X^{n-1}), 
\end{equation}
where $ \tilde{H}$ is reduced homology and $X^n/X^{n-1}$ is the identification space
obtained by collapsing the subspace $X^{n-1}$ to a point. The space  $X^n/X^{n-1}$ is
homeomorphic to a wedge of $n$-spheres, one for each $C_i^n$.  

\begin{notation}\label{not-celterm} 
Typically, for the geometric realization of an $n$-cell we shall use terminology
and the associated notation of maps on surfaces, especially when considering a 
specific case. So, we use the terms and notation for branch point, vertex, edge, face,
and polygon. On the other hand, especially when discussing the general, regular case,
it may  be simpler to let $\beta^n_i$ stand for the geometric realization of an $n$-cell
with boundary. Likewise, in a specific case, we may the use the geometrical objects 
noted above when discussing cellular boundary maps.
\end{notation}

Next, we need to define $\partial_n$ and state the cellular boundary formula 
\eqref{eq-CBF}. The differential $\partial_n$ is the composition of standard 
boundary and inclusion maps in homology:
\begin{equation}\label{dia-cellchain2} 
\xymatrix{
    \partial_n : H_n(X^n,X^{n-1}) \ar[r]^-\partial & 
    H_{n-1}(X^{n-1})\ar[r]^-j & H_{n-1}(X^{n-1},X^{n-2}).
   }
\end{equation}
Then, \emph{the cellular boundary formula} (see \cite{Ha}) is:  
\begin{equation}\label{eq-CBF}
\partial_n(\beta^n_i) = \sum_{j \in I_{n-1}} d_{i,j} \beta^{n-1}_j,
\end{equation}        
where $ d_{i,j}$ is the degree of a map to be described next. 

The cell $B^n_i$ is attached to the $n-1$ skeleton by an attaching map
\[
\alpha^{n-1}_i:\partial{B^n_i}=S^{n-1}_i \rightarrow X^{n-1}.
\]
We follow this by a quotient map 
\[
q^{n-1}_j:X^{n-1} \rightarrow S^{n-1}_j 
\]
obtained by collapsing all the spheres in $X^n/X^{n-1}$ to a point except the one
corresponding to $C^{n-1}_j$. Then  $q^{n-1}_j\circ \alpha^{n-1}_i$ is a self-map 
of the $n-1$ sphere, so we may define 
\begin{equation}\label{eq-degree}
  d_{i,j} = \mathrm{degree}(q^{n-1}_j\circ a^{n-1}_i).  
\end{equation} 

\begin{remark}\label{rk-matrep}
The formula \eqref{eq-CBF} shows that $\partial_n$ is represented by the 
integer matrix $[d_{i,j}]$.  For a regular cell complex the combined map 
$q^{n-1}_j\circ \alpha^{n-1}_i$ is a homeomorphism, resulting in $d_{i,j}=\pm 1$, 
by \eqref{eq-degree}.
\end{remark}

\section{Homology of maps on a surface}\label{sec-maphomology} 
\subsection{Computational considerations and assumptions.}
We will now impose regularity upon our surface map $(U,\mathcal{G})$. 
We need to give each $2$-cell and $1$-cell an orientation, these orientations will
determine the choices for $\beta^2_i$ and $\beta^1_j$. The orientation on the surface
defines an orientation  on each $2$-cell,  in turn telling us the sequence of edges 
and their directions in a counterclockwise circuit around the cell boundary.  
The orientations of the 1-cells may need to be picked arbitrarily since we can have 
two different induced orientations from the 2-cells that contain the edge in their 
common boundary. By regularity,  $q^1_j\circ \alpha^1_i$ is a self-homeomorphism, 
and $d_{i,j} =1$  if the counter-clockwise circuit about the boundary of the face 
labelled $i$ travels along the edge labelled $j$ in the assigned direction. 
If the direction is opposite then $d_{i,j} = -1$. 

Since our calculation of homology is going to involve linear algebra computations 
using the matrices described in  Remark \ref{rk-matrep}, we should try to make these
matrices as small as possible so that the construction of the matrices 
from the generating vector is as easy as possible. As such we construct the map $(T,\mathcal{E})$ with the following conditions in mind. 
\begin{enumerate}
  \item The lifted map is a regular $CW$-complex.
  \item The number of lifted 2-cells, 1-cells, and 0-cells are as small as practical.
  \begin{itemize}
    \item The number of 2-cells  is $|G|\times|F|$, where $|F|$ is the number 
       of faces of $(T,\mathcal{E})$.
    \item The number of 1-cells  is $|G|\times|E|$, where $|E|$ is the number 
       of edges of $(T,\mathcal{E})$.
    \item The number of 0-cells  is 
       \[ |G|\times\left(|W|+\sum_{j=1}^t\frac{1}{n_j}\right),
       \]
       where $|W|$ is the number of regular nodes in $(T,\mathcal{E})$.
  \end{itemize}
  \item The labelling and orientation of the cells, and hence the 
     determination of $d_{i,j}$, can be done efficiently,  just using 
     a generating vector. This allows for efficient coding in Sage.
\end{enumerate}

\begin{remark}\label{rk-numfaces}
We will consider two different cases for the map $(T,\mathcal{E})$,
one with a single face, ($\mathcal{E}=\mathcal{E}_1$) and the second with two faces 
 ($\mathcal{E}=\mathcal{E}_2$). The single faced case generalizes
in a nice way when the genus of $T$ is greater than 1. The two faced case 
is a better for vizualization (black and white tiles) and computation, even 
though the intermediate matrices are somewhat bigger.  
\end{remark}

We are going to make some assumptions on our vertices, edges, and faces 
of the quotient $(T,\mathcal{E})$ and the equivariant tiling $(S,\tilde{\mathcal{E}})$, 
in order to efficiently compute the cellular homology.  
The validity of the assumptions depends on choices made for
$\mathcal{E}$, a preferred face $f$ on $T$, a basepoint 
$w$ in the interior of $f$, and a generating set for $\pi_1(T^\circ,w)$.
Despite having made all these choices, we will capture all homology 
representations up to standard representation equivalence and the action of $\mathrm{Out}(G)$
on homology representations.

\begin{assumption}\label{asmpt-distinguished}
As noted previously, $G$ acts freely on the open 1-cells and 2-cells. Considering
$\beta_k^1$ and  $\beta_l^2$ as geometric cells with boundary, we may   
enumerate the cells as pairs $(g,\beta_k^1)$ and  $(g,\beta_l^2)$ 
where $\beta_k^1$ and  $\beta_l^2$ range over a set of distinguished orbit 
representatives of the 1-cells and 2-cells. The boundary maps $\partial_n$ are 
dependent on the representatives so we are going to make these assumptions.
\begin{enumerate}
  \item All the edges and vertices of $\mathcal{E}$  occur on the boundary of every face 
  in $T-\mathcal{E}$.
  
  \item Let $\beta^2$ be a fixed but arbitrary face of $S$. Then for each 1-cell $\beta^1$ (edge) 
  and 0-cell $\beta^0$ (vertex), the orbits $G\cdot\beta^1$ and $G\cdot\beta^0$ 
  have representatives that lie in the geometric boundary $\partial\beta^2$ of 
  the face $\beta^2$.  
 
  \item  For each face in  $f$ in $T-\mathcal{E}$, each branch point $z_j$ occurs 
  exactly once in the idealized boundary of $f$.  We may chose a preferred face $f$
  in $T-\mathcal{E}$ such that branch points are in the cyclic order $z_1,\ldots,z_t$, 
  when we travel counter-clockwise along the idealized boundary $\partial f$.
  
  \item For every face $\beta^2$ of $S$, a ramified vertex $v_j$, lying over $z_j$,
  occurs exactly once in $\partial\beta^2$. Moreover, for those faces $\beta^2$ 
  lying over the preferred face $f$, the vertices are in the cyclic order 
  $v_1,\ldots,v_t$, when we travel counter-clockwise along $\partial\beta^2$. 

  \item Let $(c_1,\ldots,c_t)$ be the generating vector of the action.
   Then there is a distinguished face $\beta_0^2$ such that the ramified vertices 
  $v_j \in \partial\beta_0^2$ have the following stabilizers:
\[ 
  \mathrm{Stab}_G(v_j) =\langle c_j \rangle, 
\]
  and the local action of $c_j$  at $v_j$ is a 
  counter-clockwise rotation through  $2\pi/n_j$ radians.  
\end{enumerate}   
\end{assumption}

Our assumptions may seem a little ambitious, but they are not. The construction of  $\mathcal{E}_1$ and $\mathcal{E}_2$ in Sections \ref{subsec-1face}
and \ref{subsec-2face} were specialized so that items 1 and 3 would hold.

Item 2 follows directly from item 1. To prove item 4, note that the open
faces of $S$ lift homeomorphically from the open faces of $T-\mathcal{E}$. 
In addition, the orientation on $S$, and hence upon its 2-faces may be assumed 
to be lifted from $T$ since $\pi_G$ is a holomorphic map. 

To prove item 5 we first select $w$ in the interior of our preferred face $f$ on $T$. 
We may choose a generating set for $\pi_1(T^\circ,w)$ such that the generator
$\gamma_j$ encircles $z_j$ exactly once in the counter-clockwise direction. 
These generators are described in some detail in Section \ref{subsec-2face}.
We compute all generating vectors for the various monodromies $\xi$ with respect
to this basis (see Proposition \ref{prop-xiGV equiv}).
Now pick any distinguished face $\beta_0^2$ and $\tilde{w} \in \beta_0^2$ lying over $w$. 
The lifting map $\lambda$ and hence $\epsilon$ is computed with respect to this
lifted base point. From the discussion in Section \ref{subsec-lifting} 
on lifting and base points, we see that action of $c_j =\xi(\gamma_j)$ is 
given by the lift $\lambda(\gamma_j)$ based at $\tilde{w}$. The local character
of this lift may be directly examined to prove item 5.

\subsection{The single faced case}\label{subsec-1face} 

\paragraph{\textbf{Labelling lifted cells in the single faced case.}}
Now we turn specifically to the map $(S,\widetilde{\mathcal{E}})$, 
for the cut system $\mathcal{E}=\mathcal{E}_1$  and a lifted polygon
$\mathcal{P}$ illustrated in Figure \ref{tikz-EnP}. The black nodes
 correspond to the branch points $z_j$ with the same numbering. 
The white node $z_0 \in \mathcal{E}$ is a regular point introduced to make 
the lifted map $CW$-regular. The diagram for a different number 
$t\ge 3$ of  branch points is similar, with the proper branch points 
$z_j$, $1\le j \le t$, arranged on a circle centred at  $z_0$ 
and  $j$ cyclically increasing in clockwise order. 
The edges of $\mathcal{E}$ are  the \emph{spokes} $e_j =(z_0,z_j)$, 
emanating from  $z_0$ along a radial line to the branch point $z_j$.  
Notice that in a general cut-system with a $CW$-regular lift, 
the number of lifts of a vertex lying on $\partial \mathcal{P}$ equals
the  valence of that vertex in $\mathcal{E}$. 
So, in the case at hand, each proper branch point $z_j$ lifts to a 
unique vertex on $\partial \mathcal{P}$, denoted $v_j$, and the white node $z_0$ lifts
to $t$ distinct white nodes $w_j$, $1\le j \le t$. See Figure \ref{tikz-EnP} 
for an illustration of the lifted vertices $v_j$ and $w_j$. 

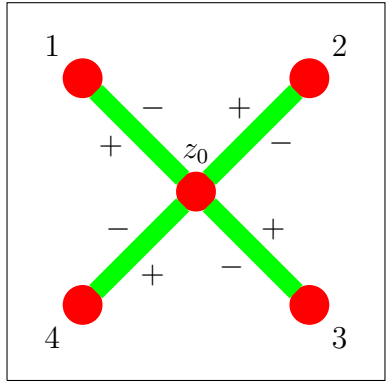
\begin{figure}[h!]
\begin{center}

\begin{tikzpicture}[scale=0.5, inner sep=1mm]

\node (A) at (3,3)   [shape=circle, fill=red, minimum size=15] [label={45:$2$}] {};
\node (B) at (-3,3)  [shape=circle, fill=red, minimum size=15] [label={135:$1$}] {};
\node (C) at (-3,-3)  [shape=circle, fill=red, minimum size=15] [label={-135:$4$}] {};
\node (D) at (3,-3) [shape=circle, fill=red, minimum size=15] [label={-45:$3$}] {};
\node (W) at (0,0)   [shape=circle, fill=red, minimum size=15] [label=$z_0$] {};

\node (E1+) at (1.5+0.4,1.5) [label={135:$+$}] {};
\node (E1-) at (1.5,1.5+0.5) [label={-45:$-$}] {};
\node (E2+) at (-1.5,1.5+0.4) [label={-135:$+$}] {};
\node (E2-) at (-1.5-0.4,1.5) [label={45:$-$}] {};
\node (E3+) at (-1.5-0.4,-1.5) [label={-45:$+$}] {};
\node (E3-) at (-1.5+0.2,-1.5-0.2) [label={135:$-$}] {};
\node (E4+) at (1.5-0.2,-1.5-0.2) [label={45:$+$}] {};
\node (E4-) at (1.5+0.2,-1.5+0.2) [label={-135:$-$}] {};

\draw [green,line width=7] (2.65,2.65)   to (0.35,0.35);
\draw [green,line width=7] (2.65,-2.65)  to (0.35,-0.35);
\draw [green,line width=7] (-2.65,2.65)  to (-0.35,0.35);
\draw [green,line width=7] (-2.65,-2.65) to (-0.35,-0.35);

\draw [thin] (-5,-5) to (5,-5) to (5,5) to (-5,5) to (-5,-5) ;

\end{tikzpicture}

\end{center}
\caption{Ribbon graph for $\mathcal{E}_1$, with orientation markings}
\label{tikz-ribbon}
\end{figure}

We follow the method and notation of \cite{BrCoIz3} to determine the 
counter-clockwise labelling of the sequence of lifted edges and nodes on 
$\partial \mathcal{P}$.  We construct a ribbon graph 
$\mathcal{R}$, which is a disc neighbourhood of $\mathcal{E}$, 
as illustrated in Figure \ref{tikz-ribbon}. 
The boundary  $\partial\mathcal{R}$, suitably oriented,
will determine the labelling and orientation of the vertices and edges in 
$\partial \mathcal{P}$. First, include in $\mathcal{R}$ a small red disc 
about each branch point and the white node, so that the system of discs 
is disjoint, and such that each red disc only meets the spokes incident 
with the node at the centre of the red disc. Next, for each spoke 
construct two line segments parallel and very close to the spoke,
on each side of the spoke. Include in $\mathcal{R}$ the area
between each pair of parallel segments including the spoke. 
These are the green ribbons in Figure \ref{tikz-ribbon}.

%%%% tikzfigure ribbon

A clockwise circuit  along $\partial \mathcal{R}$ is a counter-clockwise circuit
in the open face $T-\mathcal{E}_1$, very close to its boundary. The circuit is 
a concatenation of paths doing the following: 
\begin{itemize}
    \item start near $z_0$ and travel outward along a parallel close to the edge $e_j$,
    \item make a tight clockwise turn about the branch point $z_j$, along the boundary of the red disc, 
    \item return to $z_0$ along the opposite parallel of $e_j$.
\end{itemize}
In Figure \ref{tikz-ribbon} the direction of travel
along the edges of the green ribbons in $\partial \mathcal{R}$ is indicated by a $+$
sign (outward travel), or $-$ sign (inward travel). The lift $\widetilde{R}$ of the 
ribbon graph meets the closed polygon $\overline{\mathcal{P}}$ in a thin ribbon 
$\widetilde{R} \cap \overline{\mathcal{P}}$ that spans the space between
$\partial \mathcal{P}$ and the lifted boundary $\widetilde{\partial{R}}$.
This lifted ribbon gives us a mechanism for aligning $\partial \mathcal{P}$ 
and  $\widetilde{\partial{R}}$.

In the right hand panel of Figure \ref{tikz-EnP} the black node lying over $z_j$ 
is labeled by $v_j$.  The circuit also allows us to label the white nodes by $w_j$, so that 
in the counter-clockwise travel along $\partial\mathcal{P}$, $w_j$ lies between $v_{j-1}$ and $v_j$, 
indices mod $t$. Correspondingly, $v_j$ lies between $w_j$ and $w_{j+1}$, again indices mod $t$.  
Every edge $e_j$ in $\mathcal{E}_1$ has two lifts in $\partial\mathcal{P}$, which we denote by 
$e_j^+$ and  $e_j^-$. When travelling around the boundary of the lifted ribbon graph,
$e_j^+$ corresponds to travelling outward along the spoke $(z_0,z_j)$ and $e_j^-$ corresponds
to travelling inward along the spoke $(z_j,z_0)$. As previously described, 
in Figure \ref{tikz-ribbon} the portions corresponding to $e_j^+$ and  $e_j^-$ 
are indicated by $+$ and $-$ signs. We observe that 
\[
 v_j = e_j^+ \cap e_j^- \text{ and } w_j = e_j^+ \cap e_{j-1}^-.
\] 

Our construction shows that every lifted polygon can be given a boundary structure
and labelling that satisfies items  1-4 of Assumption \ref{asmpt-distinguished}.   
For item 5 of Assumption \ref{asmpt-distinguished}, we must now fix a distinguished 
polygon, which we still label $\mathcal{P}$. All other lifted polygons have the form 
$g\mathcal{P}$ for some $g \in G$. Having chosen $\mathcal{P}$, we must restrict the
lifting homomorphism $\lambda$ to be based at the unique point $\tilde{w} \in \mathcal{P}$ 
lying over the base point $w$. With this restriction, item 5 will hold.

Our lifted geometric objects $v_j$, $w_j$, $e_j^+$, and $e_j^-$ will now refer to the lifts of vertices 
and edges in $\widetilde{\mathcal{E}_1}$ lying in the boundary $\partial \mathcal{P}$
of the distinguished polygon. There is a \emph{side pairing transformation} $\sigma_{e_j} \in G$ such that 
$\sigma_{e_j}$ maps $e_j^+$ to $e_j^-$ in an orientation reversing manner and 
\[
e_j^- =\mathcal{P} \cap  \sigma_{e_j}\mathcal{P}.
\] 

\paragraph{\textbf{Cellular homology in the single faced case.}}
We now select our bases of  $H_n(X^n,X^{n-1})$ in direct sum format. We first need 
to enumerate the cells. As noted, the open 2-cells are $g\mathcal{P}, g \in G$. The 
open 1-cells are $ge_j^+, g \in G$, noting that $e_j^- = \sigma_{e_j}e_j^+$. In  
\cite{BrCoIz3} it was determined that  $\sigma_{e_j}=c_j^{-1}$ so that
\begin{equation}\label{eq-ejpm}
  e_j^- = c_j^{-1}e_j^+,
\end{equation}
with orientation reversed.
For  0-cells we use the quantities $w_j$ and $v_j$ to denote the lifts of vertices 
in $\mathcal{E}_1$ to the closure of $\mathcal{P}$. Consequently, we can say 
\begin{equation}\label{eq-edge def}
 e_j^+ = (w_j,v_j) \text{ and } e_j^- = (v_j,w_{j+1}),
\end{equation}
indices modulo $t$.  From equations \eqref{eq-ejpm} and   \eqref{eq-edge def} we see 
that $ w_{j+1}= c_j^{-1}w_j$ and so 
\begin{equation}\label{eq-wj}
  w_{j+1} = c_j^{-1}\cdots c_1^{-1}w_1=h_j^{-1}w_1 \text{ for } 1 <j< t, 
\end{equation} 
where $ h_j = c_1\ldots c_j$ ($h_0=1$). Thus, the 0-cells are $gw_1, g \in G$ 
and $hv_j \text{ for } h \in G/\langle c_j\rangle$. So, in direct sum format, 
we may write
\begin{eqnarray}
H_2(X^2/X^{1}) &=& \bigoplus_{g \in G} \mathbb{Z}g\mathcal{P} ,\\
H_1(X^1/X^{0}) &=& \bigoplus_{j=1}^t\bigoplus_{g \in G} \mathbb{Z} ge_j^+ ,\\
\tilde{H}_0(X^0) &=& \left( \bigoplus_{j=1}^t\bigoplus_{h \in G/\langle c_j\rangle} 
\mathbb{Z}  hv_j\right) \bigoplus \left( \bigoplus_{g \in G} \mathbb{Z}w_1\right).
\end{eqnarray}
Finally, we can write down the formulas for boundary maps
\begin{equation}\label{eq-d2}
\partial_2(g\mathcal{P}) = g\partial_2(\mathcal{P})= 
\sum_{j=1}^t ge^+_j-gc_j^{-1}e^+_j  ,
\end{equation}
and
\begin{eqnarray}\label{eq-d1}
\partial_1(ge_j^+)= g\partial_1(e_j^+) &=& gv_j-gw_j ,\\
&=& gv_j-gh_{j-1}^{-1}w_1 , \\
\partial_1(ge_j^-)= gc_j^{-1}\partial_1(e_j^+) &=& - gc_j^{-1}
v_j+gc_j^{-1}h_{j-1}^{-1}w_1 , \\
                                       &=& gh_j^{-1}w_1 -gv_j . 
\end{eqnarray} 
 
\subsection{The two faced case}\label{subsec-2face} 
We propose an alternative map  $(T,\mathcal{E}_2)$ for spherical $T$, 
it has the advantage of starting off as a $CW$-regular map. 

As in the single faced case we may assume that the branch points are $t \geq 3$ 
equally spaced points on the equator of $T$ with the same arrangement of branch points 
as in the single faced case. We choose a planar representation of $T$ by 
stereographic projection to $\mathbb{C}$ so that the south pole corresponds to the origin.  
Now we construct a cell complex on $T$ as follows. There are two $2$-cells: 
the upper and lower hemispheres. In the plane representation these 2-cells 
are the open unit disc and the interior of the complement of the unit disc, 
completed by the point at infinity.   The points on the equator determine a 
series of $t$ arcs, $e_j$, which are the $1$-cells. We chain them together to form  
$\mathcal{E}_2$. The $t$ branch points $z_j$ are the $0$-cells. 
We think of $T$ as two regular polygons with three or more sides glued together 
along corresponding edges. Over the open  $2$-cells and $1$-cells $\pi _{G}$ 
is unramified and hence in $S$ the lifts of $2$-cells and the $1$-cells are 
polygons and arcs upon which $G$ acts freely. Note that the lifts of the closed 
two cells and the closed $1$-cells map homeomorphically onto their images in $T$ 
since their vertices are all distinct.

Now we determine a specific labelling of the cells in $S$. Let $\mathcal{P}_{L}$ 
denote a specific polygon lying over the lower half sphere and 
$\tilde{e}_{1},\ldots,\tilde{e}_{t}$ the edges of $\mathcal{P}_{L}$ with a 
counter-clockwise ordering induced from the lower hemisphere. For 0-cells 
we let $v_{1},\ldots ,v_{t}$ be the vertices of $\mathcal{P}_{L}$, with 
$v_j$ lying over $z_j$.  Consequently, $\tilde{e}_{j}$ starts at $v_{j-1}$ 
and ends at $v_{j}$ (indices modulo  $t$).  The vertices and branch points 
also have a counter-clockwise ordering. There is a unique polygon 
$\mathcal{P}_{U}$ lying over the upper hemi-sphere such that $\mathcal{P}_{L}$ 
and $\mathcal{P}_{U}$  meet along the edge $\tilde{e}_{1}$.
A counter-clockwise circuit around $\mathcal{P}_{U}$ induces the opposite 
ordering on edges in $T$ and vertices and reverses orientation induced on the edges. 

The $2$-cells of $S$ have the form $g\mathcal{P}_{L}$ and $g\mathcal{P}_{U}$,
$g\in G,$ and there are $2\left\vert G\right\vert $ of them.  We call the 
$g\mathcal{P}_{L}$  lower faces and the $g\mathcal{P}_{U}$  upper faces. 
The $1$-cells have the form $g\tilde{e}_{j}$, $g\in G$ and there are 
$t\left\vert G\right\vert$ of them. We say that the edges of the 
form $g\tilde{e}_{j}$ have type $j$. The $0$-cells or vertices have the 
form $gv_{j}$ though they are not all distinct as we vary $g$. As
in the case of edges, vertices of the form $gv_{j}$ have type $j$.

Now we pick a specific base point $w=0 \in \mathcal{P}_{L}$ and generating set of $\pi_{1}(T^{\circ},w)$ to construct the generating vector from the monodromy map $\xi:\pi _{1}(T^{\circ },w)\twoheadrightarrow G.$ 
Let $\gamma _{j}$ be the loop that starts at $w$, follows an imaginary spoke to
the branch point corresponding to $z_{j}$, makes as a small counter-clockwise
loop around the branch point and heads back along the spoke to $w$.
 Our generating vector $(c_{1},\ldots ,c_{t})$ is given by $c_{j}=\xi (\gamma_{j})$. 
With these choices, the element $c_{j}$ fixes $v_{j}$ and permutes 
the faces with vertex corner $v_{i}$ by a counterclockwise rotation
through angle $2\pi /n_j$. The lower and upper faces are placed in an
alternating fashion around $v_{i}$. 

For homology computations it is important for us to know the labels of the $t$ edges
that form the boundary of $\mathcal{P}_{U}$. By simple transitivity, for each edge 
$\tilde{e}_{j}$ of $\mathcal{P}_L$ there is a unique $k_{j}\in G$ such that 
the lower face $k_j\mathcal{P}_L$ meets $\mathcal{P}_{U}$ along 
$k_j\tilde{e}_{j}$. We may determine $k_j$ by path lifting. 
Let $\gamma$ be a loop in $T^\circ$ starting at $w$, crossing the edge 
$e_1$ transversally into the upper hemisphere, and then crossing back into 
the lower hemisphere through the edge $e_j$  back to $w$. 
The path crosses edges only as specified and lies entirely in $T^\circ$. 
Consider the lifted path $\tilde{\gamma}$ that starts at $\tilde{w} \in \mathcal{P}_L$,
lying over $w$, crosses $\tilde{e}_1$ into  $\mathcal{P}_{U}$, crosses 
$k_j\tilde{e}_j$ back into $k_j\mathcal{P}_L$ and ends up at $\tilde{\gamma} (1) =\xi(\gamma)\tilde{w}=k_j\tilde{w}$.  It is not hard to show that 
$\gamma =\gamma_1\cdots\gamma_{j-1}$ so that 
\[
k_j=\xi(\gamma)=\xi(\gamma_1)\cdots\xi(\gamma_{j-1})=c_1\cdots c_{j-1}=h_{j-1},
\]
as previously defined.

Here are the boundary maps
\begin{eqnarray}
  \partial_2(g\mathcal{P}_L) &=&  \sum_{j=1}^t g\tilde{e}_j, \\
  \partial_2(g\mathcal{P}_U) &=&  -\sum_{j=1}^t gk_j\tilde{e}_j, \\
   \partial_1(g\tilde{e}_j)  &=&  gv_j-gv_{j-1} . 
\end{eqnarray}

\section{The homology calculation algorithm}
Our homology calculation and the $G$ action $g\rightarrow \rho(g)$ upon it 
is now a linear algebra problem to be implemented as a computer calculation in Sage. 
By definition 
$H_{1}(S;\mathbb{Z})=\mathop{ker} (\partial _{1})/\mathop{im}(\partial _{2})$. 
The chain group $H_1(X^1/X^0)$ is a free $\mathbb{Z}$ module of rank $t|G|$ with 
a $\mathbb{Z}$  basis indexed by $t$ copies of $G$. We must identify a basis for 
$\mathop{ker}(\partial_1)$ and then a basis for a complementary subspace for 
$\mathop{im}(\partial _{2})$. We propose to do this by Gaussian elimination on 
large matrices $M_{1}$, $M_{2}$ representing the operators $\partial _{1}$  and 
$\partial_{2}$ respectively. 

To organize our vectors and matrices, we first make an ordered list of the elements of 
$G$: $g_1=1, g_2,\ldots, g_{|G|}$. Then our ordering of components of vectors corresponding to 1-cells
will be based an ordering of the pairs $(g_i,e_j)$ lexicographic on the the pair $(j, i)$, where $e_j$ is an edge in $\mathcal{E}$. The same remarks apply to 2-cells and 0-cells but some care is needed for the 0-cells
as transversals of  $G/\langle c_j\rangle$ will be required. 

\subsection{Some linear algebra}\label{subsec-linaalg}
A key step of our process will be selecting a basis from a spanning set of a subspace. 
Here is a basic lemma concerning such a selection which is true for vector 
spaces over fields. It follows from a careful examination of Gaussian elimination 
as taught in a first course in linear algebra.
\begin{lemma}\label{lem-findbasis}
  Let 
\[
M=
  \left[
  \begin{array}{ccc}
    C_1 & \cdots & C_n 
  \end{array}
  \right]
\]
be a matrix written as a row of column vectors. Let $M^\prime$ be the matrix $M$ after 
performing Gaussian elimination, using row operations. Let $P\subseteq \{1,\ldots,n\}$ 
be the set of indices of the pivot columns of  $M^\prime$. Then, the set of vectors
$\{C_i: i \in P\}$ is a basis of the column space  of  $M$. Moreover, $m \in P$ if 
and only if $C_m$ is not in the linear span of the columns that precede it.
\end{lemma}  

Our next lemma will aid us in computing the $G$ action, but first we set up 
the notation and a bit of background for the lemma. Let $V\subset W$ be 
two subspaces of an ambient vector space, and let
\begin{equation}\label{eq-RCD}
 R=
  \left[
  \begin{array}{cccccc}
    C_1 & \cdots & C_m&D_1& \cdots&D_n  
  \end{array}
  \right]
  =
  \left[ 
 \begin{array}{cc}
    C&D  
  \end{array}
  \right]
\end{equation}
be a matrix with independent columns such that $C_1,\dots,C_m$ is a basis of 
$V$ and $W$ is the column space of $R$.  We can set up local coordinates on $V$ 
and $W/V$ by writing for $Y\in W$: 
\begin{equation}\label{eq-localcoord}
 Y=CY_C+DY_D.
\end{equation}
where $Y_C,Y_D$ are unique column vectors of the length $m$ and $n$ respectively.  
It is not hard to show that the maps defined by
\begin{equation}\label{eq-localcoordiso}
Y_C\rightarrow CY_C \text{ and }  Y_D\rightarrow DY_D+V
\end{equation}
are isomorphisms onto $V$ and $W/V$, respectively.

Now, Let $L$ be a matrix that induces a linear transformation $Y\rightarrow LY$ 
in the ambient space.  Also denoting the transformation by $L$, we assume 
that $L(V) \subseteq V$ and  $L(W) \subseteq W$,  and hence $L$ induces a 
linear transformation $L_Q$ on $W/V$. Next we write 
\[
 LCY_C=CY^\prime_C \text{ and }  LDY_D= CY^{\prime\prime}_C + DY^\prime_D.
\]
There is no $D$ term in the first equation since $V$ is $L$-invariant.
The maps $Y_C\rightarrow Y^\prime_C$ and $Y_D\rightarrow  Y^\prime_D$ are 
the induced maps on $V$ and $W/V$, in local coordinates. 
The first statement is clear, for the second we note that at the coset level
\[
L_Q(DY_D +V) =  CY^{\prime\prime}_C + DY^\prime_D + V = DY^\prime_D + V.
\]
We want to find matrices $L_C$, $L_D$ of appropriate size such that 
\begin{equation}\label{eq-localtrans1}
Y^\prime_C=L_CY_C \text{ and }  Y^\prime_D=L_DY_D.
\end{equation}
Inspired by the method of finding the normal equations for least squares, 
we need to find compatible left inverses for $C$ and $D$. Indeed, as
$R= \left[ \begin{array}{cc} C&D \end{array} \right]$  has independent columns, 
it has a left inverse: 
\begin{equation}\label{eq-leftinv}
\left[ \begin{array}{c} A\\ B \end{array} \right]
\left[ \begin{array}{cc} C&D \end{array} \right]
=\left[ \begin{array}{cc} AC&AD \\ BC&BD \end{array} \right]
=\left[ \begin{array}{cc} I_m&0 \\ 0&I_n \end{array} \right],
\end{equation}
for suitable matrices $A,B$.
Then,
\[
 ALCY_C=ACY^\prime_C=Y^\prime_C
\]
and
\[
 BLDY_D= BCY^{\prime\prime}_C + BDY^\prime_D = Y^\prime_D.
\]
We select
\begin{equation}\label{eq-localtrans2}
  L_C = ALC \text{ and } L_D = BLD
\end{equation} 

The preceding discussion proves the following lemma.
\begin{lemma}\label{lem-findtrans}
Let notation be as above. Then:
\begin{enumerate}
  \item The images of $D_1\ldots,D_n$ in $W/V$ form a basis for $W/V$.
  \item There are local coordinate vectors $Y_C$ and $Y_D$ for $V$ and $W/V$ such that    the maps in \eqref{eq-localcoordiso}  are isomorphisms onto $V$ and  
     $W/V$.  
  \item Let $L$ be a matrix inducing a linear transformation $L: Y\rightarrow LY$ 
     on the ambient space, that satisfies $L(V)=V$ and $L(W)=W$. 
     Let  $A,B$ be the components of a left inverse of $R$ as defined 
     in \eqref{eq-leftinv}. Then in the local coordinates, $Y_C$ and $Y_D$, 
     the induced linear transformation of $L$ on $V$ and $W/V$ are given by 
     equations \eqref{eq-localtrans1} and \eqref{eq-localtrans2} 
\end{enumerate} 
\end{lemma}

\subsection{Steps for computing homology representation}
Here now are the steps for finding the homology representation.

\begin{enumerate}
  \item Set up vector coordinates for the chain group $H_1(X^1/X^0)$ 
     as suggested in the opening paragraphs of this section.
  \item Calculate the matrices $M_1$ and $M_2$ of the differentials $\partial_1$ 
     and $\partial_2$, also as described in the opening paragraphs of this section. 
  \item Using the rows of $M_2$, and Lemma \ref{lem-findbasis}, find a basis for 
     $\mathop{im}(\partial_2)$. Fill in the $C$ portion of the matrix $R$ 
     in \eqref{eq-RCD} with this basis. 
  \item Find a basis for the null space of $M_1$.   
  \item Make a temporary matrix  
     $R^\prime= \left[ \begin{array}{cc} C&D^\prime \end{array} \right]$, 
     where $D^\prime$ is filled in with the basis from the previous step.
  \item Using Gaussian elimination on $R^\prime$ as described in Lemma 
     \ref{lem-findbasis}, prune out the columns of $D^\prime$ that are 
     dependent on the columns of $C$, resulting in 
     $R= \left[ \begin{array}{cc} C&D \end{array} \right]$.
  \item Find a left inverse of $R$ as in equation \eqref{eq-leftinv}. 
     The left inverse can be found using the standard method of finding 
     right inverses using Gaussian elimination of an augmented matrix applied to $R^T$.
  \item Select $g \in G$ for which a homology transformation matrix 
     $\rho(g)$ is desired. Create the matrix 
     \[ 
     R^g= \left[ \begin{array}{cc} C^g&D^g \end{array} \right]
     \] 
  by applying the transformation induced by $g$ on $H_1(X^1/X^0)$
  to all the columns of $R$. The matrix $L$ of Lemma \ref{lem-findtrans} 
  need not be created,  we need only create the modified columns in 
  $C^g$ and $D^g$.
  \item Use Lemma \ref{lem-findtrans} to find the desired matrix.
\end{enumerate} 

\begin{remark}\label{rk-gencoeff}
According to Hatcher \cite{Ha}, cellular homology with coefficients works 
as described above with some attention paid to the cellular boundary 
formula \eqref{eq-CBF}. In that formula, $\beta^{n-1}_j$ lies in the 
coefficient group and  $d_{i,j} \beta^{n-1}_j$ is simply the $\mathbb{Z}$ 
module action of multiplying a coefficient by an integer. The remainder of this 
section works without change when the coefficient group is a field, 
in particular $\mathbb{Z}_2$ coefficients.
\end{remark}

\section{The intersection form and matrix}
%Mention \S\ref{subsec-homrep}

The intersection form $a\smallfrown b$ on a surface $S$ is a skew-symmetric $R$-bilinear form on $H_1(S;R)$, where $R$ 
is any unital commutative ring. The form is invariant under the action of  $G$ on $S$. Once a basis $a_1,\ldots,a_{2\sigma}$ of $H_1(S;R)$ is chosen we may construct the \emph{intersection matrix}
\begin{equation*}
I_S=
\left[ 
\begin{array}{c}
a_i\smallfrown a_j
\end{array}
\right].
\end{equation*}
Suppose $X = (x_1,\ldots,x_{2\sigma})$ and $Y=(y_1,\ldots,y_{2\sigma})$ are two $R$-vectors and
$a_X = x_1a_1+\cdots+x_{2\sigma}a_{2\sigma}$ and $a_Y =y_1a_1+\cdots+y_{2\sigma}a_{2\sigma}$ are the corresponding
homology classes. Then 
\[
a_X\smallfrown a_Y = X^\top I_SY.
\]
Through a change of basis the intersection matrix may be brought to canonical form $\begin{pmatrix} 0 & I \\ -I & 0 \end{pmatrix}$, whereby the matrices of the homology representation become symplectic, taking values in $Sp_{2\sigma}(R)$. As this shall be necessary for the application to theta characteristics in \S\ref{sec-theta}, we shall briefly describe how this may be implemented, restricting to the two faced case of \S\ref{subsec-2face}. 

It may be shown that if $a^\top$ and  $b^\top$ are the Poincar\'{e} duals of 
$a,b$, respectively then 
\[ 
a\smallfrown b = (a^\top\smallsmile b^\top) \smallfrown[S],
\]
where $\smallsmile$ denotes the cup product in cohomology and $[S]$ is the fundamental 
homology $2$-class of $S$. An explicit formula for the cup product in simplicial cohomology is given in \cite[\S3.2]{Ha}, and we may use this this in our implementation provided we refine the $CW$-complex on $S$ to a simplicial complex. In the case $t=3$ there is no change to the complex. In the $t>3$ case this refinement may be achieved by adding additional edges through the upper and lower faces, see Appendix \ref{app-refined-CW}. This does increase the dimension of the matrices which must be used in the computation of the homology basis. Refining the complex by adding additional edges (and so also additional faces) has the effect of modifying the differentials $\partial_{1,2}$ from those defined earlier, but the remaining method for computing the homology basis and the action of $G$ on it carries through. 

By using the cup product formula we avoid the problem present in other works of having to find representatives of homology classes, and avoid difficulties in correctly signing intersections. As computational implementation of the cup product in simplicial cohomology is well known, for example see \cite{Ga}, we shall not discuss the details further here. 

\section{Examples and computation}

\subsection{Hyperelliptic example}

As a demonstration of how the procedure works, we shall apply the algorithm to $S$ determined by $G = \langle \iota \, | \, \iota^2=1\rangle$ acting with signature $(0; 2, 2, 2, 2)$ and generating vector $(\iota, \iota, \iota, \iota)$. This is the action of the hyperelliptic involution on a genus-1 surface.

\begin{remark}
    Note this surface $S$ is not hyperbolic, but the group action considered is still faithful and the method still applies. 
\end{remark}

$M_2$ is the $t|G| \times 2|G|$ matrix (in the case of the hyperelliptic involution $|G|=2$ and $t= 2\sigma+2$) which acts on column vectors representing 2-cells ordered as 
$$
g_1 \mathcal{P}_L, \dots, g_{|G|}\mathcal{P}_L, g_1 \mathcal{P}_U, \dots, g_{|G|}\mathcal{P}_U, 
$$ 
by left multiplication. The output column vectors correspond to 1-cells ordered as 
\[
g_1 \tilde{e}_1, \dots, g_{|G|}\tilde{e}_1, g_1 \tilde{e}_2, \dots, g_{|G|}\tilde{e}_t. 
\]
As such, introducing the $|G|\times|G|$ matrices $\rho_i = \rho_R(k_i)$ the images of $k_i$ under the right-multiplication representation of $G$ on itself, we find 
\[
M_2 = \begin{pmatrix}
    I_{|G|} & -\rho_1 \\ \vdots & \vdots \\ I_{|G|} & -\rho_t
\end{pmatrix}, 
\]
where $I_n$ is the $n\times n$ identity matrix. 

Likewise, if we introduce the $\frac{|G|}{n_i} \times |G|$ matrices $P_i$ which act by left multiplication on column vectors of length $|G|$ to project elements $g \in G$ to their left coset $g\langle c_i \rangle$, then we have $M_1$ is the matrix 
\[
M_1 = \begin{pmatrix}
P_1 & -P_1 & 0 & \dots & 0  \\ 0 & P_2 & -P_2 & \dots& 0 \\ \vdots & \vdots & \ddots & \ddots & 0 \\ 0 & \vdots & 0 & P_{t-1} & -P_{t-1} \\ -P_t & 0 & \dots  & 0 & P_t
\end{pmatrix}. 
\]

In the case at hand we see 
\[
\rho_i = \left \lbrace \begin{array}{cc}
    I_2, & i \text{ odd}, \\
    \begin{pmatrix} 0 & 1 \\ 1 & 0 \end{pmatrix}, & i \text{ even}, 
\end{array} \right. \quad \text{and} \quad P_i = \begin{pmatrix} 1 & \dots & 1 \end{pmatrix}. 
\]

Following the procedure laid out one ends up finding that 
\[
D = \begin{pmatrix}
    1 & 0 \\ 0 & 0 \\ 0 & 1 \\ 1 & -1 \\ 0 & 0 \\ 1 & 0 \\ 0 & 0 \\ 1 & 0 
\end{pmatrix}, \quad B = \begin{pmatrix}
    1 & 0 & 0 & 0 & -1 & 0 & 0 & 0 \\ 0 & 0 & 1 & 0 & 0 & 0 & -1 & 0
\end{pmatrix}. 
\]
Moreover matrix $L$ corresponding to the action of $g$ is 
\[
L = \bigoplus_{i=1}^{t} \rho_L(g),
\]
where $\rho_L$ is the left-multiplication representation of $G$ on itself. As such $\rho_L(\iota) = \begin{pmatrix} 0 & 1 \\ 1 & 0 \end{pmatrix}$, and one may verify $L_D = -I_2$ as is to be expected for the homology action of the hyperelliptic involution.
% \[ 
% L = \mathrm{diag}\left( \begin{pmatrix} 0 & 1 \\ 1 & 0 \end{pmatrix},\begin{pmatrix} 0 & 1 \\ 1 & 0 \end{pmatrix},\begin{pmatrix} 0 & 1 \\ 1 & 0 \end{pmatrix},\begin{pmatrix} 0 & 1 \\ 1 & 0 \end{pmatrix} \right). 
% \]

\subsection{Implementation comments}

The algorithm is implemented in Sage using the complex of the two faced case, and so being valid when $\tau=0$. It is available from \url{https://github.com/DisneyHogg/homology_representation/}. 

Some comments about the implementation are appropriate. 
\begin{itemize}
    \item GAP is used to get the permutation representation of $G$ acting on the cosets $G/\langle c_i \rangle$. 
    \item The matrices $M_{1}, M_2$ are implemented as sparse matrices in an endeavor to keep the algorithm as memory efficient as possible when $|G|$ gets large. This proves to be a sensible choice; our implementation is able to compute the representation of $\mathrm{Aut}(S)$ on $H_1(S; \mathbb{Z}_2)$ for the modular curves $S= X(13)$, $X(17)$, which existing algorithms are unable to due to memory issues (see \S\ref{sec-theta}). 
    \item It is convenient when computing the intersection form to work directly with cohomology rather than homology, but the material effect of this is just to work with row vectors from the beginning rather than column vectors and later transposing.  
\end{itemize}

\subsection{Performance comparison}

In order to evaluate the performance of our new algorithm, we wish to compare to existing methodologies (especially memory efficiency), namely that from \cite{BeRoRo}. For both algorithms we found the time taken to compute the representation on $H_1(S; \mathbb{Z})$ of the $G$ action for all 1897 topologically inequivalent generating vectors acting in genera $2 \leq \sigma \leq 12$, the latter data generated using \cite{Behn2023}. These runtimes will naturally be hardware dependent; computations were run on an Intel Core i5-8350U CPU at 1.70GHz. At this level we see that our new algorithm is on average $1.2\times$ faster. 

Moreover, one key motivation of this paper was to have a methodology which performed better when computing the representation on $H_1(S; \mathbb{Z}_2)$, as this is required for computing orbits of theta characteristics (see the later discussion in \S\ref{sec-theta}). As such we recomputed the runtime in this scenario and again compare the performance, seeing now that the computation is on average $2\times$ faster. 

Somewhat more importantly, as the genus of the curve increases our new method appears to have better memory efficiency, rendering new ranges of genera amenable to computation. We ran both algorithms to compute the homology representation of $C_p$ acting with signature $(0; p, p, p)$ for increasing $p$ an odd prime, which acts on a surface of genus $\sigma = \frac{p-1}{2}$. The method of \cite{BeRoRo} was unable to compute the representation for $p \geq 79$ due to overflow errors, but our new algorithm continued to work and answer in less than one second up to $p=331$.

\subsection{Theta characteristics}\label{sec-theta}

As previously mentioned, an important application for us of this method will be computing the action of $G$ on $H_1(S; \mathbb{Z}_2)$, because of its role in computing the orbits of theta characteristics \cite{Kallel2010, Braden2025}; these are line bundles $L \to S$ such that $L^{\otimes 2} \cong K_S$ which may be identified with elements of $H_1(S; \mathbb{Z}_2)$, such that the induced action of $\mathrm{Aut}(S)$ on theta characteristics becomes an affine action of the homology representation on $H_1(S; \mathbb{Z}_2)$. 

In particular we consider the modular curves $X(p)$, $p \geq 7$ an odd prime. These are curves of genus $\sigma = \frac{1}{24}(p+2)(p - 3)(p - 5)$ that have an action of $PSL(2, p)$ with signature $(0; 2,3,p)$. The results developed in \cite{Braden2025} yield that these curves should have a unique invariant characteristic, and we seek to verify that numerically where possible. While the code of \cite{BeRoRo} works in the cases $p=7, 11$, it does not complete for larger $p$. Using our new method we can check $p=13, 17$, verifying the theoretical results. 

Moreover, with the homology representation in hand we can compute the parity of the invariant characteristic, which at present is unable to be determined by theory alone in most cases. We find that the unique invariant characteristic always has even parity for $p = 7, 11, 13, 17$: we conjectured that this was always the case, a result which has subsequently been proven \cite{Disneyhogg20XX}.

%%%%%%%%%%%%%%%%%%%%%%%%%%%%%%%%%%%%%%%%%%%%
%%%%%%%%%%%%%%%%%%%%%%%%%%%%%%%%%%%%%%%%%%%%
\appendix

%%%%%%%%%%%%%%%%%%%%%%%%%%%%%%%%%%%%%%%%%%%%
%%%%%%%%%%%%%%%%%%%%%%%%%%%%%%%%%%%%%%%%%%%%

\section{Details of 2-face simplicial complex}\label{app-refined-CW}

Here we give further details of how the $CW$-complex constructed in the 2-face case with $t>3$ is refined to give a simplicial complex in order to use the cup product to compute the intersection form. 

From the map $(T, \mathcal{E}_2)$ we construct a map $(T, \mathcal{E}_2^\prime)$ by adding arcs from $z_1$ to $z_i$ ($i=3, \dots, t-1$) in each hemisphere, denoted $l_i$, $u_i$. Lifting this as before the 1-cells on $S$ are $g\tilde{e}_j$ ($j=1, \dots, t$) and $g \tilde{l}_i$, $g \tilde{u}_i$ ($i=3, \dots t-1$). The new 1-cells have differentials 
\[
\partial_1(g \tilde{l}_i) = g v_i - gv_1, \quad \partial_1(g \tilde{u}_i) = g k_i v_i - g k_1 v_1 .
\]

The 2-cells on $S$ are now $g L_i$, $g U_i$, ($i=3, \dots, t$) which are lifts of faces on the lower and upper hemisphere respectively with differentials
\[
\partial_2(g L_i) = g \tilde{e}_i - g\tilde{l}_i + g\tilde{l}_{i-1}, \quad \partial_2(g U_i) = g k_i \tilde{e}_i - g \tilde{u}_i + g \tilde{u}_{i-1} ,
\]
where we have defined $g \tilde{l}_2 := g \tilde{e}_2$, $g \tilde{l}_t := - g \tilde{e}_1$, $g \tilde{u}_2 = g k_2 \tilde{e}_2$, $g \tilde{u}_t = - g k_1 \tilde{e}_1$. 
We have the interpretation in our mind that 
\[
g \mathcal{P}_L = \sum_{i=3}^t g L_i, \quad g \mathcal{P}_U = - \sum_{i=3}^t g U_i. 
\]

%%%%%%%%%%%%%%%%%%%%%%%%%%%%%%%%%%%%%%%%%%%%
%%%%%%%%%%%%%%%%%%%%%%%%%%%%%%%%%%%%%%%%%%%%
%%%%%%%%%%%%%%%%%%%%%%%%%%%%%%%%%%%%%%%%%%%%
%%%%%%%%%%%%%%%%%%%%%%%%%%%%%%%%%%%%%%%%%%%%

\bibliography{library}
% \addcontentsline{toc}{chapter}{Bibliography}
% Choose a reference style from https://verbosus.com/bibtex-style-examples.html
\bibliographystyle{abbrv}

\end{document}